\documentclass[12pt]{article}
\usepackage{amssymb}
\usepackage{latexsym,bm}
\usepackage{graphicx}
\usepackage{amsmath}
\usepackage{mathrsfs}
\usepackage{mathrsfs,amscd,amssymb,amsthm,amsmath,bm,graphicx,psfrag,subfigure,url,xcolor}
\usepackage{tikz}

\usepackage{epstopdf}

\date{}
\newcounter{mathitem}
\newenvironment{mathitem}
{\begin{list}{{$(\roman{mathitem})$}}{
\setcounter{mathitem}{0}
\usecounter{mathitem}
\setlength{\topsep}{0pt plus 2pt minus 0pt}
\setlength{\parskip}{0pt plus 2pt minus 0pt}
\setlength{\partopsep}{0pt plus 2pt minus 0pt}
\setlength{\parsep}{0pt plus 2pt minus 0pt}
\setlength{\leftmargin}{35pt}
\setlength{\itemsep}{0pt plus 2pt minus 0pt}}}
{\end{list}}
	
\begin{document}
\title{Weakly pancyclic vertices in dense nonbipartite graphs}
\author{\hskip -10mm Yurui Tang and  Xingzhi Zhan\thanks{Corresponding author.}}
\maketitle
\footnotetext[1]{Department of Mathematics,  Key Laboratory of MEA (Ministry of Education)
 and Shanghai Key Laboratory of PMMP, East China Normal University, Shanghai 200241, China}
\footnotetext[2]{E-mail addresses: {tyr2290@163.com} (Y. Tang),  {zhan@math.ecnu.edu.cn (X. Zhan).}}

\begin{abstract}
Let $G$ be a graph of girth $g$ and circumference $c.$ A vertex $v$ of $G$ is called weakly pancyclic if $v$ lies on an $\ell$-cycle for every integer $\ell$ with $g\le \ell\le c.$ We prove that if $G$ is a nonbipartite graph of order $n\ge 5$ and size at least $\left\lfloor(n-1)^2/4\right\rfloor+2,$ then $G$ contains three weakly pancyclic vertices, with one exception. This strengthens a result of Brandt from 1997. We also pose a related problem.
\end{abstract}


\vskip 8mm

\section{Introduction}

We consider finite simple graphs and use standard terminology and notations from [2] and [7]. The {\it order} of a graph is its number of vertices, and the
{\it size} is its number of edges.  A $k$-cycle is a cycle of length $k.$ A graph $G$ of order $n$ is called {\it pancyclic} if for every integer $k$ with $3\le k\le n,$ $G$ contains a $k$-cycle. See the book [4] for this topic. A vertex (edge) $v$ of a graph $G$ of order $n$ is called {\it pancyclic} if for every integer $k$ with $3\le k\le n,$ $v$ lies on
a $k$-cycle. There are natural more general concepts.

{\bf Notation 1.} We use $g(G)$ and $c(G)$ to denote the girth and circumference of a graph $G,$ respectively.

{\bf Definition 1.} A graph $G$ is called {\it weakly pancyclic} if for every integer $k$ with $g(G)\le k\le c(G),$ $G$ contains a $k$-cycle. A vertex (edge) $u$ of a
graph $G$ is called {\it weakly pancyclic} if for every integer $k$ with $g(G)\le k\le c(G),$ $u$ lies on a $k$-cycle.

Clearly, if a graph contains a weakly pancyclic vertex, then it is weakly pancyclic. The converse is false in general. The first graph in Figure 1 is pancyclic with no pancyclic
vertex while the second graph in Figure 1 has girth $4$ and circumference $7$ and it is weakly pancyclic with no weakly pancyclic vertex.
\begin{figure}[h]
\centering
\includegraphics[width=0.8\textwidth]{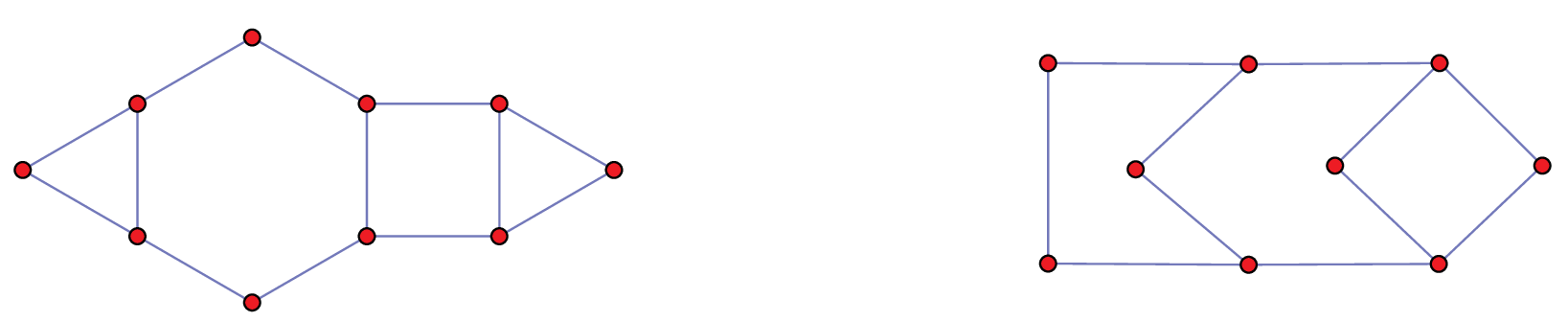}
\caption{Weakly pancyclic graphs with no weakly pancyclic vertex}
\end{figure}

The following result was conjectured by Erd\H{o}s and proved by H${\mathrm {\ddot{a}}}$ggkvist, Faudree and Schelp [5].
				
{\bf Theorem 1} [5]. {\it If $G$ is a hamiltonian nonbipartite graph of order $n$ and size at least $\left\lfloor(n-1)^2/4\right\rfloor+2,$ then $G$ is pancyclic.}		
			
Theorem 1 has been generalized by Brandt [3] as follows.
				
{\bf Theorem 2} [3]. {\it If $G$ is a nonbipartite graph of order $n$  and size at least $\left\lfloor(n-1)^2/4\right\rfloor+2,$ then $G$ is weakly pancyclic.}

It is known [3, Lemma 2] that every nonbipartite graph $G$ of order $n$ and size at least $\left\lfloor(n-1)^2/4\right\rfloor+2$ contains a triangle; i.e., $g(G)=3.$
		
We denote by $K_{s,t}$ the complete bipartite graph whose partite sets have cardinality $s$ and $t,$ respectively, and by $C_k$ the $k$-cycle.	For graphs we will use equality up to isomorphism, so $G=H$ means that $G$ and $H$ are isomorphic.

{\bf Notation 2.} For an integer $n\ge 5,$ we denote by $BT(n)$ the graph obtained by identifying an edge of $K_{\lfloor(n-1)/2\rfloor,\lceil(n-1)/2\rceil}$ with an edge of $C_3$.

The graphs $BT(8)$ and $BT(9)$ are depicted in Figure 2.

\begin{figure}[h]
\centering
\includegraphics[width=0.85\textwidth]{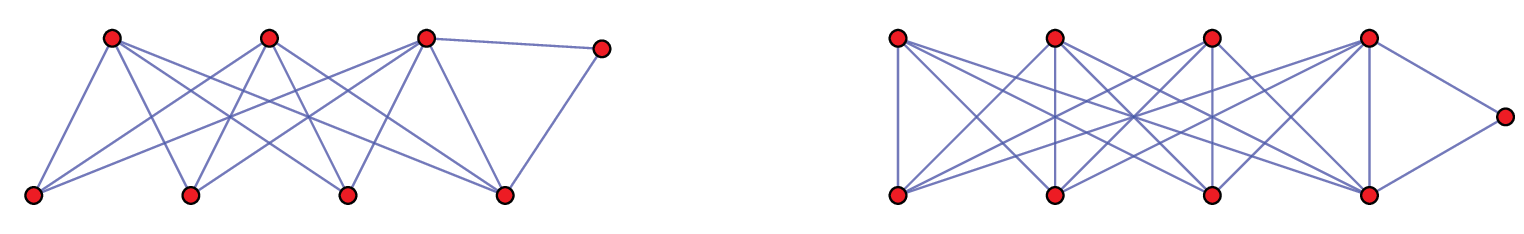}
\caption{The graphs $BT(8)$ and $BT(9)$}
\end{figure}

Observe that if $n$ is odd then $BT(n)$ contains exactly two pancyclic vertices while if $n$ is even then $BT(n)$ has circumference $n-1$ and it contains exactly two weakly pancyclic vertices.

In this paper, we strengthen Theorem 2 by proving the following result.

{\bf Theorem 3.} {\it Let $G$ be a nonbipartite graph of order $n\ge 5$ and size at least $\left\lfloor(n-1)^2/4\right\rfloor+2$ other than $BT(n).$  Then $G$ contains three weakly pancyclic vertices.}

The next result shows that the number three in Theorem 3 is sharp.

{\bf Theorem 4.} {\it For every integer $n\ge 6,$ there exists a nonbipartite graph of order $n$ and size $\left\lfloor(n-1)^2/4\right\rfloor+2$ that contains
exactly three weakly pancyclic vertices.}
				
We denote by $V(G),$ $|G|,$ $e(G)$ and $\delta(G)$ the vertex set, order, size and minimum degree of a graph $G,$ respectively.
For $v\in V(G)$, we denote by $N_G(v)$ and ${\rm deg}_G(v)$ the neighborhood and degree of $v$ in $G$, respectively.
If the graph is clear from the context, sometimes we omit the subscript $G.$
For two distinct vertices $x$ and $y$ in $G$, an {\it $(x, y)$-path}  is a path whose endpoints are $x$ and $y$.

In Section 2 we give proofs of the two main results Theorems 3 and 4, and in Section 3 we pose a related unsolved problem.

\section{Proofs}

Let $v$ be a vertex of a graph $G$ of order $n.$ We call $v$ a {\it small vertex} of $G$ if
$$
{\rm deg}_G(v)\le \lfloor (n-1)/2\rfloor=\left\lfloor (n-1)^2/4\right\rfloor-\left\lfloor (n-2)^2/4\right\rfloor;
$$ otherwise, $v$ is a {\it big vertex} of $G.$ Thus, if $n=2k+1$ or $n=2k+2,$ then $v$ is a small vertex if and only if ${\rm deg}_G(v)\le k.$
We will need the following three lemmas.
				
{\bf Lemma 5} [3, Lemma 1]. {\it If $C$ is a longest cycle of a nonhamiltonian graph $G,$ then $G-V(C)$ contains a small vertex of $G.$}

{\bf Lemma 6} [3, Lemma 2]. {\it Every nonbipartite graph of order $n$ and size at least $\left\lfloor (n-1)^2/4\right\rfloor+2$  contains a triangle.}
				
A bipartite graph of order $n$ is called {\it balanced} if its two partite sets have cardinality $\lfloor n/2\rfloor$ and $\lceil n/2\rceil,$ respectively.

{\bf Lemma 7}. {\it Let $H$ be a balanced bipartite graph of order $2k$ with  $k\ge 2$ and size at least $k^2-k+2$ with partite sets $V_1$ and $V_2.$ Then
\begin{mathitem}
\item for any vertices $u\in V_1$, $w\in V_2$ and any integer $\ell$ with $2\le \ell\le k,$  $H$ contains a $(u,w)$-path of length $2\ell-1$;
\item for any two distinct vertices $u,v\in V_1$ and any integer $\ell$ with $1\le \ell\le k-1,$  $H$ contains a $(u,v)$-path of length $2\ell.$
\end{mathitem}}

{\bf Proof}.  Part (i) can be found in [1, Lemma 6]. Next we deduce (ii) from (i).
Since $e(H)\ge k^2-k+2$, we have $\delta(H)\ge 2$ and $N(u)\cap N(v)\ne \emptyset.$ Thus (ii) holds with $\ell=1$.
We choose a vertex $x\in N(v)$ such that ${\rm deg}(x)$ is as small as possible. Note that
$$
e(H)-e(H-\{v,x\})={\rm deg}(v)+{\rm deg}(x)-1.\eqno (1)
$$

Case 1. $e(H-\{v,x\})\ge (k-1)^2-(k-1)+2$.

Since ${\rm deg}(v)\ge 2$, there exists a vertex $y\in N(v)\setminus \{x\}.$ By (i), for every integer $\ell$ with $2\le \ell\le k-1,$ there exists a $(u,y)$-path $P$ of length $2\ell-1$ in $H-\{v,x\}.$ Then $P\cup yv$ is a $(u,v)$-path of length $2\ell$.

Case 2. $e(H-\{v,x\})\le (k-1)^2-(k-1)+1$.

Now $e(H)-e(H-\{v,x\})\ge 2k-1,$ which, together with (1), implies that ${\rm deg}(v)={\rm deg}(x)=k$.  By the choice of $x$, it follows that every neighbor of $v$ has degree $k$. Then $H$ is a complete  balanced bipartite graph, and clearly part (ii) holds. \hfill $\Box$
				
Let $C=v_1v_2\cdots v_mv_1$ be a cycle. For a positive integer $p,$ we set
\begin{align*}
\text{
$v_i^+=v_{i+1}$,~~
$v_i^-=v_{i-1}$,~~
$v_i^{+p}=v_{i+p}$~~and~~
$v_i^{-p}=v_{i-p}$,
}\end{align*}
where the subscripts are to be read modulo $m.$ For a subset $S\subseteq V(C)$, we denote
\begin{align*}
S^{+p}=\{v_i^{+p}|\, v_i\in S\}\text{~~and~~}
S^{-p}=\{v_i^{-p}|\, v_i\in S\}.
\end{align*}
If $i<j$, we denote by $v_i\overrightarrow{C}v_j=v_iv_{i+1}\cdots v_{j}$ and $v_i\overleftarrow{C}v_j=v_iv_{i-1}\cdots v_{j}$ the two $(v_i, v_j)$-paths on $C.$
Similar notations will be used for paths.				
				
{\bf Proof of Theorem 3.}
We first prove the case when $G$ is hamiltonian; i.e., the following
				
{\bf Proposition 8.} {\it Let $G$ be a hamiltonian nonbipartite graph of order $n\ge 5$ and size at least $\left\lfloor(n-1)^2/4\right\rfloor+2$ other than $BT(n).$  Then $G$ contains three pancyclic vertices.}
				
We use induction on $n$ to prove Proposition 8. If $n=5$, or $n=6$, it is easy to verify that the conclusion of Proposition 8 holds. Next let $n\ge  7$ and assume that Proposition 8 holds for all graphs of order less than $n.$
	
{\bf Claim 1.} If $G$ contains a small vertex $x$ such that $G-x$ is hamiltonian, then Proposition 8 holds.				
				
{\bf Proof of Claim 1.} Let  $C=v_1v_2\cdots v_{n-1}v_1$ be a Hamilton cycle of $G-x$ and denote $v_n=x.$
By the definition of a small vertex, we have $e(G-v_n)\ge \left\lfloor (n-2)^2/4\right\rfloor+2.$ If $G-v_n$ is neither a bipartite graph nor $BT(n-1)$, by the induction hypothesis, $G-v_n$ contains three pancyclic vertices $u_1, u_2, u_3.$ Since $G$ is hamiltonian, $u_1, u_2, u_3$ are also pancyclic vertices of $G.$
				
Suppose that $G-v_n$ is a bipartite graph. Since $G-v_n$ is hamiltonian and bipartite, $n$ is odd. Let $n=2k+1$ where $k\ge3,$ and let $V_1$ and $V_2$ be the two partite sets of $G-v_n$ with $|V_1|=k$ and $|V_2|=k$. By Lemma 6, $G$ contains a triangle. Since $G-v_n$ is bipartite, $v_n$ lies on every triangle. Suppose $v_nyzv_n$ is a triangle where $y\in V_1$ and $z\in V_2$.
Note that $e(G-v_n)\ge k^2-k+2$. By Lemma 7(i), for every integer $\ell$ with $2\le \ell\le k,$ there exists a $(y,z)$-path $P_1$ of length $2\ell-1$ in $G-v_n.$ Then $P_1\cup yz$ is a cycle of length $2\ell$ and $P_1\cup yv_nz$ is a cycle of length $2\ell+1$. Hence, $y$ and $z$ are pancyclic vertices of $G$. Since $G\ne BT(n)$, we have ${\rm deg}_G(v_n)\ge 3$. Let $w\in N_G(v_n)\setminus\{y,z\}$. Without loss of generality, suppose $w\in V_1$. By Lemma  7(ii), for every integer $\ell$ with $1\le \ell\le k-1,$
there exists a $(y,w)$-path $P_2$ of length $2\ell $ in $G-v_n.$ Then we find a cycle $P_2\cup yv_nw$ of length $2\ell+2$, $1\le\ell\le k-1$, that is, $v_n$ lies on cycles of all even lengths. We have already shown that $v_n$ lies on cycles of all odd lengths. Thus $v_n$ is a pancyclic vertex of $G,$ and $G$ contains three pancyclic vertices $v_n$, $y$ and $z.$
				
Next suppose $G-v_n=BT(n-1)$. Recall that $BT(p)$ is hamiltonian if and only if $p$ is odd.
Now $n$ is even, since $G-v_n$ is hamiltonian. Let $n=2k+2$ where $k\ge 3.$ Without loss of generality, suppose that $v_{n-1}$ is the vertex of $G-v_n$ that has degree $2,$ and $N_{G-v_n}(v_{n-1})=\{v_{n-2},v_1\}$.  Let $V_1$ and $V_2$ be the partite sets of the bipartite graph $G-\{v_n,v_{n-1}\}=K_{k,k}$ such that $v_1\in V_1$ and $v_{n-2}\in V_2.$
Note that $v_1$ and $v_{n-2}$ are pancyclic vertices of $G-v_{n}$, which implies that $v_1$ and $v_{n-2}$ are also two pancyclic vertices of $G,$ since $G$ is hamiltonian.

We use the notation $w_1\leftrightarrow w_2$ to mean that two vertices $w_1$ and $w_2$ are adjacent, and use the notation $w_1\nleftrightarrow w_2$ to mean that
$w_1$ and $w_2$ are nonadjacent. With $n=2k+2,$ $e(G)\ge \left\lfloor(n-1)^2/4\right\rfloor+2=k^2+k+2.$ Since $e(G-v_n)=k^2+2,$ we have ${\rm deg}_G(v_n)\ge k\ge 3.$

 We first assume that $v_n\leftrightarrow v_{n-1}$. Observe that $v_{n-1}$ lies on cycles of all lengths $\ell$ in $G-v_n,$ where $\ell$ is an odd integer with $3\le\ell \le 2k+1.$
If $v_n\leftrightarrow v_1$ or $v_n\leftrightarrow v_{n-2}$, we obtain an $(\ell+1)$-cycle that contains $v_{n-1}$ by replacing the edge $v_{n-1}v_1$ or $v_{n-1}v_{n-2}$ with $v_{n-1}v_nv_1$ or $v_{n-1}v_nv_{n-2}$ in an $\ell $-cycle. Then  $G$ contains three pancyclic vertices $v_{n-1}$, $v_{n-2}$ and $v_1$. Now suppose
$v_n\nleftrightarrow v_1$ and $v_n\nleftrightarrow v_{n-2}.$ Since ${\rm deg}_G(v_n)\ge 3$, there exists a vertex $w\in N_G(v_n)\setminus\{v_{n-1}\}$. Without loss of generality, suppose $w\in V_2$. Since $v_1$ and $w$ are adjacent and $G-\{v_n,v_{n-1}\}=K_{k,k},$ there exists a $(v_1,w)$-path $P_3$ of length $2\ell-1$ in $G-\{v_{n-1},v_n\}$ for
$1\le \ell \le k$. Then $P_3\cup wv_nv_{n-1}v_1$ is a cycle of length $2\ell+2$. Clearly $v_{n-1}$ lies on cycles of all odd lengths.
Thus, $v_{n-1}$ is a pancyclic vertex of $G.$ Consequently $G$ contains three pancyclic vertices $v_{n-1}$, $v_{n-2}$ and $v_1$.

Next we assume that $v_n\nleftrightarrow v_{n-1}$. Since every hamiltonian graph is tough and $G$ is hamiltonian, we have $N_G(v_n)\cap V_1\ne\emptyset$ and $N_G(v_n)\cap V_2\ne\emptyset.$ Since ${\rm deg}_G(v_n)\ge 3,$ there exist $w_1\in N_G(v_n)\cap V_1$ and $w_2\in N_G(v_n)\cap V_2$ such that $\{w_1, w_2\}\neq \{v_1, v_{n-2}\}.$ Then $w_1\leftrightarrow w_2.$ Observe that $BT(n-1)$ is a spanning subgraph of $G-v_{n-1}.$ Using a similar argument as in the above proof that $v_1$ and $v_{n-2}$ are pancyclic vertices of $G,$ we deduce that $w_1$ and $w_2$ are pancyclic vertices of $G.$ Since $\{w_1, w_2\}\neq \{v_1, v_{n-2}\},$ the set $\{w_1, w_2\}\cup\{v_1, v_{n-2}\}$ contains at least three pancyclic vertices of $G.$
This completes the proof of Claim 1.		
			
{\bf Claim 2.} If $G$ contains a big vertex $y$ such that $G-y$ is hamiltonian, then $y$ is a pancyclic vertex of $G.$
				
{\bf Proof of Claim 2.} Let $n=2k+1$ if $n$ is odd and let $n=2k+2$ if $n$ is even where $k\ge 3.$ Let $D$ be a Hamilton cycle of $G-y.$ Since $y$ is a big vertex of $G$,
 ${\rm deg}_D(y)={\rm deg}_{G}(y)\ge k+1.$ For any integer $\ell$ with $3\le \ell\le n,$ we have $|N_D(y)|=|[N_D(y)]^{+(\ell-2)}|\ge k+1$ and $|D|\le 2k+1.$ Hence there exists a vertex
  $u\in N_D(y)\cap [N_D(y)]^{+(\ell-2)}.$ Then $u^{-(\ell-2)}\overrightarrow{D}u\cup uyu^{-(\ell-2)}$ is an $\ell$-cycle of $G.$ This completes the proof of Claim 2.

By Theorem 1, $G$ contains an $(n-1)$-cycle $C=v_1v_2\cdots v_{n-1}v_1.$  Let $V(G)\setminus V(C)=\{v_n\}.$ If $v_n$ is a small vertex, then by Claim 1, Proposition 8 holds.
Next we treat the case when $v_n$ is a big vertex of $G$. By Claim 2, $v_n$ is a pancyclic vertex of $G$.

Case 1. $n$ is odd.

Let $n=2k+1$ with $k\ge3.$ Then ${\rm deg}_C(v_n)\ge k+1$. Since $|C|=2k$, we have
$$
|N_{C}(v_n)\cap [N_{C}(v_n)]^{+2}|\ge 2,
$$
which implies that there exist two distinct vertices $w_1$ and $w_2$ such that $\{w_1^-,w_2^-,w_1^+,w_2^+\}\subseteq N_C(v_n)$. Consider the two distinct $(n-1)$-cycles
$$v_nw_1^+\overrightarrow{C}w_1^-v_n~~\text{and}~~v_nw_2^+\overrightarrow{C}w_2^-v_n.$$
If $w_1$ or $w_2$ is a small vertex of $G$, by Claim 1, Proposition 8 holds. So we assume that both $w_1$ and $w_2 $ are big vertices of $G$. By Claim 2, $w_1$ and $w_2$ are two pancyclic vertices of $G.$ Thus $G$ contains three pancyclic vertices $v_n, w_1, w_2.$
				
Case 2. $n$ is even.

Let $n=2k+2$ with $k\ge3.$  Then ${\rm deg}_C(v_n)\ge k+1$. Similar to the case of $n=2k+1$, if $|N_{C}(v_n)\cap [N_{C}(v_n)]^{+2}|\ge 2,$ then Proposition 8 holds, so assume that $|N_{C}(v_n)\cap [N_{C}(v_n)]^{+2}|\le 1.$ Since ${\rm deg}_C(v_n)\ge k+1$, we deduce that $|N_{C}(v_n)\cap [N_{C}(v_n)]^{+2}|= 1$ and ${\rm deg}_C(v_n)= k+1$.
Suppose $w_3$ is the vertex of $C$ such that $w_3^+\leftrightarrow v_n $ and $w_3^-\leftrightarrow v_n$. Consider the $(n-1)$-cycle $C'=w_3^-v_nw_3^+\overrightarrow{C}w_3^-$. If $w_3$ is a small vertex of $G$, by Claim 1, Proposition 8 holds. Next we assume that $w_3$ is a big vertex of $G.$ Since $G-w_3$ is hamiltonian, by Claim 2, $w_3$ is a pancyclic vertex of $G.$ If $|N_{C'}(w_3)\cap [N_{C'}(w_3)]^{+2}|\ge 2,$ then Proposition 8 holds as above, so assume that $|N_{C'}(w_3)\cap [N_{C'}(w_3)]^{+2}|$ $\le 1.$
Since ${\rm deg}_{C^{\prime}}(w_3)\ge k+1$, we deduce that $|N_{C^{\prime}}(w_3)\cap [N_{C^{\prime}}(w_3)]^{+2}|= 1$ and ${\rm deg}_{C^{\prime}}(w_3)= k+1.$

Subcase 2.1. $v_n\nleftrightarrow w_3.$

Let $P=w_3^+\overrightarrow{C}w_3^-.$ We have  ${\rm deg}_P(v_n)={\rm deg}_C(v_n)=k+1$.
Suppose $N_{P}(v_n)=\{f_1,f_2,\dots,f_{k+1}\}$ where $f_1,f_2,\dots,f_{k+1}$ appear on $C$ in order, $f_1=w_3^+$ and $f_{k+1}=w_3^-$.
Since $|N_{C}(v_n)\cap [N_{C}(v_n)]^{+2}|= 1,$ we have $\{f_1^{+2},f_2^{+2},\dots,f_{k-1}^{+2}\}\cap N_P(v_n)=\emptyset.$ Since $|N_P(v_n)|=k+1$ and $|C|=2k+1,$ we have
$V(C)=\{f_1^{+2},f_2^{+2},\dots,f_{k-1}^{+2}\}\cup N_C(v_n)\cup\{w_3\},$  implying that $k+1$ is even and
$$
N_C(v_n)=\{w_3^{+},w_3^{+2},w_3^{+5},w_3^{+6},\dots, w_3^{-2},w_3^{-}\}=\mathop{\cup}\limits_{i=1}^{(k+1)/2}\{w_3^{+(4i-3)}, w_3^{+(4i-2)}\}.
$$
Interchanging the roles of $v_n$ and $w_3,$ using the above argument we obtain $N_G(w_3)=N_G(v_n).$ In particular,
$w_3\leftrightarrow w_3^{+2}.$ Observe that $G- w_3^{+}$ contains a Hamilton cycle $w_3w_3^-v_nw_3^{-2}\overleftarrow{C} w_3^{+2}w_3$. If $w_3^{+}$ is a small vertex of $G,$ by Claim 1, Proposition 8 holds; otherwise by Claim 2, $w_3^{+}$ is a pancyclic vertex of $G$ and thus $G$ contains three pancyclic vertices $v_n, w_3, w_3^{+}.$

Subcase 2.2. $v_n\leftrightarrow w_3.$

Since $|N_{C}(v_n)\cap [N_{C}(v_n)]^{+2}|= 1$ and $\{w_3^-, w_3, w_3^{+}\}\subseteq N_G(v_n),$ we deduce that
$$
\{w_3^{-2}, w_3^{-3}, w_3^{+2}, w_3^{+3}\}\cap N_G(v_n)=\emptyset,
$$
implying $k\ge 4,$ since  ${\rm deg}_C(v_n)=k+1.$ Consider the path $Q=w_3^{+4}\overrightarrow{C}w_3^{-4}.$ Note that $|Q|=2k-6$ and ${\rm deg}_Q(v_n)=k-2.$
Using a similar argument as in Subcase 2.1, we deduce that $k$ is even and
\begin{align*}
N_C(v_n)&=\{w_3^-, w_3, w_3^{+}\}\cup \{w_3^{+4},w_3^{+5},w_3^{+8},w_3^{+9},\dots, w_3^{-5},w_3^{-4}\}\\
        &=\{w_3^-, w_3, w_3^{+}\}\cup\mathop{\cup}\limits_{i=1}^{(k/2)-1}\{w_3^{+(4i)}, w_3^{+(4i+1)}\}.
\end{align*}
Now we verify that $w_3^{+}$ is a pancyclic vertex of $G.$ First $w_3^{+}$ lies on the triangle $v_n w_3 w_3^{+}v_n.$ Next for an integer $\ell$ with $4\le\ell\le 2k-1,$
if $\ell=4i,$ $w_3^{+}$ lies on the $\ell$-cycle $v_nw_3^{-}\overrightarrow{C}w_3^{+(4i-3)}v_n;$  if $\ell=4i+1,$ $w_3^{+}$ lies on the $\ell$-cycle $v_nw_3^{+}\overrightarrow{C}w_3^{+4i}v_n;$ if $\ell=4i+2,$ $w_3^{+}$ lies on the $\ell$-cycle $v_nw_3^{+}\overrightarrow{C}w_3^{+(4i+1)}v_n;$ if $\ell=4i+3,$ $w_3^{+}$ lies on the $\ell$-cycle $v_nw_3\overrightarrow{C}w_3^{+(4i+1)}v_n.$ Finally
$w_3^{+}$ lies on the $(2k)$-cycle $v_nw_3^{-}\overrightarrow{C}w_3^{-4}v_n,$ on the $(2k+1)$-cycle $C,$ and on a $(2k+2)$-cycle, since $G$ is hamiltonian.

Thus $G$ contains three pancyclic vertices $v_n, w_3, w_3^{+}.$  This completes the proof of Proposition 8.							
		
Now we use induction on $n$ to prove Theorem 3. If $n=5$ or $n=6$, it is easy to verify that Theorem 3 holds. Next let $n\ge  7$ and assume that Theorem 3 holds for all graphs of order less than $n.$ If $G$ is hamiltonian, then Proposition 8 shows that the conclusion of Theorem 3 holds. Next we assume that $G$ is nonhamiltonian.
				
Let $R$ be a longest cycle of $G$. By Lemma 5, there exists a small vertex $x$ in $G-V(R)$. Denote $M=G-x.$
Then $e(M)\ge \left\lfloor (n-2)^2/4\right\rfloor+2.$
				
{\bf Claim 3.} $M$ is nonbipartite.
				
{\bf Proof of Claim 3.}
To the contrary, suppose $M$ is a bipartite graph with partite sets $V_1,V_2$ where $|V_1|\ge |V_2|.$ Let $n_1=|V_1|$ and $n_2=|V_2|$. Clearly, $n_1+n_2=n-1$
and $|R|\le 2n_2.$ By Lemma 6, $G$ contains a triangle. Since $M$ is bipartite, $x$ lies on every triangle.
Suppose $xyzx$ is a triangle where $y\in V_1$ and $z\in V_2$. Let $H$ be a graph obtained from $M$ by deleting $n_1-n_2$ vertices in $V_1\setminus\{y\}$.
Then $H$ is a balanced bipartite graph of order $2n_2.$  Note that for any vertex $w\in V(M)\setminus V(H)$, ${\rm deg}_{M}(w)\le n_2$.
Recall that $e(M)\ge \left\lfloor (n-2)^2/4\right\rfloor+2.$ Then
\begin{align*}
e(H)&\ge e(M)-(n_1-n_2)n_2
\\&\ge \left\lfloor \frac{(n-2)^2}{4}\right\rfloor+2-(n_1-n_2)n_2
\\&= \left\lfloor \frac{(n_1+n_2+1)^2-4(n_1+n_2+1)+4}{4}\right\rfloor+2-n_1n_2+n_2^2
\\&= \left\lfloor \frac{(n_1-n_2-1)^2}{4}\right\rfloor+2+n_2^2-n_2
\\&\ge n_2^2-n_2+2.
\end{align*}
By Lemma 7(i), there exists a $(y,z)$-path $P$ of length $2n_2-1$ in $H.$ Thus we find a cycle $P\cup yxz$ of length $2n_2+1>2n_2\ge |R|$, a contradiction.
This completes the proof of Claim 3.

{\bf Claim 4.} $M\ne BT(n-1).$

{\bf Proof of Claim 4.} To the contrary, assume $M=BT(n-1).$ If $n$ is even, let $n=2k+2$ where $k\ge 3$. Let $p$ be the vertex of  $M$ of degree $2,$
and let $V_1,V_2$ be the partite sets of $M-p=K_{k,k}.$ Since $e(G)\ge \left\lfloor(n-1)^2/4\right\rfloor+2=k^2+k+2$ and $e(M)=k^2+2$, we have ${\rm deg}_G(x)\ge k$. Note that $BT(n-1)$ is hamiltonian if $n$ is even.  Using the condition
$c(G)=c(M)=n-1,$ it is easy to verify that $N_G(x)=V_1$ or 	$N_G(x)=V_2,$ which implies that $G=BT(n),$ contradicting the assumption of Theorem 3. If $n$ is odd,  let $n=2k+1$ where $k\ge 3$. Since $e(G)\ge k^2+2$ and $e(M)=k^2-k+2,$ we have ${\rm deg}_G(x)\ge k$. It follows that  $c(G)>n-2=c(M)=c(G),$ a contradiction.
This completes the proof of Claim 4.

By Lemma 6, $M$ has girth $3.$ Since $g(M)=g(G)$ and $c(M)=c(G)$, every weakly pancyclic vertex of $M$ is also a weakly pancyclic vertex of $G.$

Using Claims 3, 4 and the fact that $e(M)\ge \left\lfloor (n-2)^2/4\right\rfloor+2,$ and applying the induction hypothesis to $M,$ we deduce that $M$ contains three weakly pancyclic vertices which are also weakly pancyclic vertices of $G.$ This completes the proof of Theorem 3. \hfill $\Box$
				
{\bf Proof of Theorem 4.} Denote $H=BT(n-1)$ and let $x$ be the vertex of $H$ with degree $2$ that lies in the triangle. Let $V_1, V_2$ be the partite sets of the balanced bipartite graph $H-x$ with $|V_1|=\lfloor (n-2)/2\rfloor$ and $|V_2|=\lceil (n-2)/2\rceil.$ Let $G_n$ be the graph obtained from $H$ by adding a new vertex $y$ and joining $y$ to $x$ and to
every vertex in $V_2\setminus N_H(x).$ The graphs $G_9$ and $G_{10}$ are depicted in Figure 3. Observe that $G_n$ is obtained from $BT(n)$ by deleting one edge and adding one edge.
\begin{figure}[h]
\centering
\includegraphics[width=0.7\textwidth]{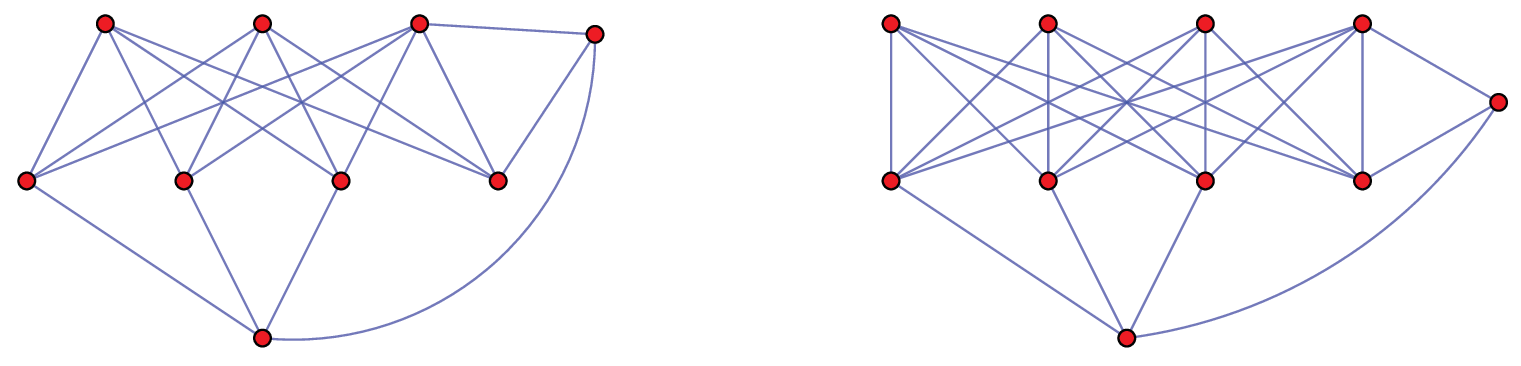}
\caption{The graphs $G_9$ and $G_{10}$}
\end{figure}

Now we show that $G_n$ is what we want. First, $G_n$ has order $n$ and size $\lfloor(n-1)^2/4\rfloor+2,$ and it is nonbipartite since it contains a triangle.
Clearly $G_n\ne BT(n),$ since the unique triangle of $G_n$ contains no vertex of degree $2.$ By Theorem 3, $G_n$ contains three weakly pancyclic vertices. On the other hand,
since $G_n-x$ is a bipartite graph, every vertex of $G_n$ other than the three in the triangle does not lie on any triangle. Hence the three vertices of the triangle
are the only weakly pancyclic vertices of $G_n.$ This completes the proof of Theorem 4.\hfill $\Box$

\section{Unsolved problems}

We first mention a related problem which can be found in [6].
				
{\bf Conjecture 1} (X. Zhan, July 2025) [6]. If $G$ is a hamiltonian nonbipartite graph of order $n\ge 7$ and size at least $\left\lfloor(n-1)^2/4\right\rfloor+2$ other than $BT(n)$ when $n$ is odd, then $G$ contains a pancyclic edge.

Every endpoint of a pancyclic edge is a pancyclic vertex. Thus Conjecture 1, if true, would imply that every hamiltonian nonbipartite graph of order $n\ge 7$ and size at least $\left\lfloor(n-1)^2/4\right\rfloor+2$ other than $BT(n)$ when $n$ is odd contains two pancyclic vertices.

Finally we pose a new problem.

{\bf Problem 2.} For a positive integer $n,$ let $f(n)$ denote the smallest integer $k$ such that every nonbipartite graph of order $n$ and size at least $k$ contains a weakly pancyclic vertex. Determine $f(n).$

Denote $b(n)=\left\lfloor(n-1)^2/4\right\rfloor+2.$ A computer search gives the following information: $f(n)=b(n)$ for $n=6,7,8;$ $f(9)=17<18=b(9);$ $f(10)=20<22=b(10).$

\vskip 5mm
{\bf Acknowledgement.} This research  was supported by the NSFC grant 12271170 and Science and Technology Commission of Shanghai Municipality
 grant 22DZ2229014.

\section*{\normalsize Declaration}
				
\noindent\textbf{Conflict~of~interest}
The authors declare that they have no known competing financial interests or personal relationships that could have appeared to influence the work reported in this paper.
		
\noindent\textbf{Data~availability}
No data was used for the research described in the article.

\end{document}